\documentclass[12pt,reqno,oneside]{amsart}
\usepackage{amsmath}
\usepackage{amssymb}
\usepackage{amsthm}
\usepackage{graphics}
\usepackage{eucal}
\usepackage{mathrsfs}
\usepackage{bm}
\usepackage[pdftex,bookmarks,colorlinks,breaklinks]{hyperref}  
\newtheorem{theorem}{Theorem}[section]

\newtheorem{proposition}{Proposition}[section]

\newtheorem{corollary}{Corollary}[section]

\numberwithin{equation}{section}

\newcommand{\bbR}{\mathbb{R}}

\newcommand{\bbZ}{\mathbb{Z}}
\newcommand{\bt}{\mathbf{t}}
\newcommand{\bq}{\mathbf{q}}
\newcommand{\bx}{\mathbf{x}}

\newcommand{\cL}{\mathcal{L}}
\newcommand{\bp}{\mathbf{p}}

\newcommand{\bl}{\mathbf{l}}

\newcommand{\La}{\Lambda}

\newcommand{\f}{\mathbf{f}}

\newcommand{\by}{\mathbf{y}}

\newcommand{\ba}{\mathbf{a}}
\newcommand{\bb}{\mathbf{b}}
\newcommand{\bX}{\mathbf{X}}
\newcommand{\bT}{\mathbf{T}}
\newcommand{\bS}{\mathbf{S}}
\newcommand{\ttt}{\mathbf{t}}
\newcommand{\bs}{\mathbf{s}}
\newcommand{\bll}{\mathbf{l}}
\newcommand{\bepsilon}{\boldsymbol{\varepsilon}}
\newcommand{\btheta}{\boldsymbol{\theta}}
\newcommand{\bbN}{\mathbb{N}}
\newcommand{\cM}{\mathcal{M}}
\newcommand{\cE}{\mathcal{E}}
\newcommand{\cA}{\mathcal{A}}
 
\begin{document}
\title{On Diophantine transference principles}
\begin{abstract}
We provide an extension of the transference results of Beresnevich and Velani connecting homogeneous and inhomogeneous Diophantine approximation on manifolds and provide bounds for inhomogeneous Diophantine exponents of affine subspaces and their nondegenerate submanifolds.
\end{abstract}

\subjclass[2000]{11J83, 11K60} \keywords{Diophantine approximation on manifolds, inhomogeneous Diophantine approximation, transference principles}

\author{Anish Ghosh}
\address{\textbf{Anish Ghosh} \\
School of Mathematics,
Tata Institute of Fundamental Research, Mumbai, India 400005}
\email{ghosh@math.tifr.res.in}

\author{Antoine Marnat}
\address{
\textbf{Antoine Marnat}\\
Institute of Analysis and Number Theory\\
Technische Universit\"at Graz, 8010 Graz, Austria}
\email{marnat@math.tugraz.at}

\thanks {Ghosh is supported by an ISF-UGC grant. Marnat is supported by the Austrian Science Fund (FWF), Project F5510-N26}

\maketitle

\section{Introduction}

In \cite{BV1}, V. Beresnevich and S. Velani proved beautiful transference principles which allow one to move between homogeneous and inhomogeneous Diophantine approximation on manifolds, and more generally, a class of measures introduced in their work, called \emph{contracting measures}. In a companion paper \cite{BV2}, they give a simplified version of their proof for the case of simultaneous Diophantine approximation on manifolds. We begin with this setup and then move on to a more general setting. For a vector $\bx \in \bbR^n$, let
\begin{equation}\label{defexp0}
w_{0}(\bx):= \sup\{w~:~ \|q\bx\| < |q|^{-w} \text{ for infinitely many } q \in \bbN\}
\end{equation}
and
\begin{equation}\label{defexpn-1}
w_{n-1}(\bx):= \sup\{w~:~ \|\bq\cdot\bx\| < \|\bq\|^{-w} \text{ for infinitely many } \bq \in \bbZ^{n} \backslash \{0\}\}.
\end{equation}

The exponent $w_0(\bx)$ is referred to as the simultaneous Diophantine exponent and $w_{n-1}(\bx)$ as the dual Diophantine exponent. Here and henceforth, we will use $\|x\|$ to denote the fractional part of a real number $x$, and $\|\bq\|$ to denote the supremum norm of a vector $\bq \in \bbR^n$, i.e. vectors and matrices will be denoted in boldface and for $\bq = (q_1, \dots, q_n),$
$$ \|\bq\| = \max_{1 \leq i \leq n}|q_i| $$
It is a consequence of Dirichlet's pigeon hole principle that $w_0(\bx) \geq 1/n$ and that $w_{n-1}(\bx) \geq n$ for all $\bx \in \bbR^n$. On the other hand, it is a consequence of the Borel-Cantelli lemma, that $w_0(\bx) = 1/n$ and $w_{n-1}(\bx) = n$ for Lebesgue almost every $\bx \in \bbR^n$. Similarly, in the context of inhomogeneous Diophantine approximation, one has  two analogous exponents. Since we will primarily be concerned with the simultaneous exponent, we only define its inhomogeneous counterpart. For $\btheta \in \bbR^n$, 
\begin{equation}\label{definexp0}
w_{0}(\bx, \btheta):= \sup\{w~:~ \|q\bx + \btheta\| < |q|^{-w} \text{ for infinitely many } q \in \bbN\}.
\end{equation}

Diophantine approximation on manifolds is concerned with the question of whether typical Diophantine properties in $\bbR^n$, i.e. those which are generic for Lebesgue measure, are inherited by proper submanifolds. A manifold $\cM$ is called \emph{extremal} if almost every point on $\cM$ is not very well approximable, or equivalently, if $w_0(\bx) = 1/n$ and $w_{n-1}(\bx) = n$ for almost every $\bx \in \cM$. If $\cM = \{f(\bx)~|~\bx \in U\}$ is a $d$ dimensional sub manifold of $\bbR^n$, where $U$ is an open subset of $\bbR^d$ and $f := (f_1, \dots, f_n)$ is a $C^m$ imbedding of $U$ into $\bbR^n$ and $l \leq m$, we say that $y = f(\bx)$ is an $l$-nondegenerate point of $\cM$ if the space $\bbR^n$ is spanned by partial derivatives of $f$ at $\bx$ of order up to $l$. The manifold $\cM$ will be called nondegenerate if $f(\bx)$ is nondegenerate for almost every $\bx \in U$. It was a long standing conjecture of Sprind\v{z}uk, that smooth nodegenerate manifolds are extremal. This was proved by Kleinbock and Margulis \cite{KM} in a landmark paper. Sprind\v{z}uk's formulation of his conjecture was slightly less general, the above notion of nondegeneracy is due to Kleinbock and Margulis. We refer the reader to \cite{KM} for all the details. It is natural to enquire about inhomogeneous versions of Sprind\v{z}uk's conjecture and other homogeneous results in Diophantine approximation. A manifold $\cM$ is called simultaneously inhomogeneously extremal if for every $\btheta \in \bbR^n$, 
$$w_{0}(\bx, \btheta) = \frac{1}{n} \text{ for almost every } \bx \in \cM.$$
In \cite{BV1}, Beresnevich and Velani proved the following striking theorem using a transference principle.
\begin{theorem}\label{BV}
A smooth manifold $\cM$ is extremal if and only if it is simultaneously inhomogeneously extremal.
\end{theorem}
One direction of the above Theorem is clear of course, the other, namely extremal implies simultaneously inhomogeneously extremal is the main surprise. In fact their results are much more general, and this framework is developed in the next section. The main content of \cite{BV1, BV2} is to provide an upper bound on the Diophantine exponent, namely the corresponding lower bound is provided using a transference inequality of Bugeaud and Laurent \cite{BL}. 
\begin{theorem}[Beresnevich, Velani, 2010]  
Let $\cM$ be a differentiable submanifold of $\bbR^n$. If $\cM$ is extremal, then for every $\btheta \in \bbR^n$ we have that
\begin{equation}
w_{0}(\bx, \btheta) \leq \frac{1}{n} \text{ for almost all } \bx \in \cM.
\end{equation}
\end{theorem}

For a Borel measure $\mu$ define its Diophantine exponent by 
\begin{equation}
w_0(\mu) := \sup\{v ~|~ \mu\{\bx~|~ w_0(\bx) > v\} > 0\}.
\end{equation}

\noindent The definition only depends on the measure class of $\mu$. We can similarly define the inhomogeneous exponent of a measure as follows: for $\btheta \in \bbR^n$
\begin{equation}
w_0(\mu, \btheta) := \sup\{v ~|~ \mu\{\bx~|~ w_0(\bx, \btheta) > v\} > 0\}.
\end{equation}

If $\mathcal{M}$ is a smooth submanifold of $\mathbb{R}^n$ parametrised by a smooth map $f$, then set the Diophantine exponent $w_0(\mathcal{M})$ to be equal to $w_0(f_{*}\lambda)$ where $f_{*}\lambda$ is the push forward of Lebesgue measure $\lambda$ by $f$. Then a manifold $\cM$ is extremal when $w_0(\mathcal{M}) = 1/n$ and simultaneously inhomogeneously extremal when $w_0(\mu, \btheta) = 1/n$ for all $\btheta \in \bbR^n$. The purpose of this note is to demonstrate that the method of Beresnevich and Velani can in fact be used to relate homogeneous and inhomogeneous Diophantine approximation on manifolds even when the exponent is $v \neq n$. Examples of non-extremal manifolds are given by affine subspaces and their nondegenerate manifolds. The study of Diophantine approximation of affine subspaces goes back to Schmidt and Sprind\v{zuk} and has seen significant developments recently, we refer the reader to the survey \cite{G-survey} for details. A systematic study of extremality and Diophantine exponents for affine subspaces was initiated by Kleinbock  in two beautiful papers \cite{Kleinbock-extremal, Kleinbock-exponent}. In particular, in \cite{Kleinbock-exponent}, the following result about Diophantine exponents of affine subspaces and their nondegenerate submanifolds was proved.

\begin{theorem}[Kleinbock \cite{Kleinbock-exponent}]
If $\cL$ is an affine subspace of $\bbR^{n}$ and $\cM$ is a nondegenerate submanifold in $\cL$, then
\begin{equation}
\omega_{n-1}(\cM) = \omega_{n-1}(\cL) = \inf \{ \omega_{n-1}(\bx) \mid \bx \in \cL \} = \inf\{ \omega_{n-1}(\bx) \mid \bx \in \cM \} 
\end{equation}
\end{theorem}

Furthermore, if $\cL\subset\bbR^{n}$ is a hyperplane parametrized by
\begin{equation}\label{defL}
(x_{1}, x_{2}, \cdots , x_{n-1}) \to (a_{1}x_{1} + \cdots + a_{n-1}x_{n-1}, x_{1}, \ldots , x_{n-1}),
\end{equation}

then a formula for the exponent was obtained by Kleinbock \cite{Kleinbock-exponent}.

\begin{theorem}\label{KHypn-1}
Let $\cL$ be a hyperplane defined by $\ba := (a_{1}, \ldots , a_{n-1}) \in \bbR^{n-1}$ as in \eqref{defL}. Then we have
\begin{equation}
\omega_{n-1}(\cL) = \max\left(n,\omega_{0}(\ba)\right).
\end{equation}
\end{theorem}

Here, the notion of nondegeneracy in an affine subspace is a natural extension of the definition above. Namely  if $\cL$ is an affine subspace of $\bbR^n$, $U$ is an open subset of $\bbR^d$ and $f : U \to \bbR^n$ is a differentiable map, then $f$ is said to be nondegenerate in $\cL$ at $x_0 \in U$ if $f(U) \subset \cL$ and the span of all the partial derivatives of $f$ up to some order is the linear part of $\cL$.  if  In \cite{Zhang1}, Y. Zhang provided the simultaneous analogue of Kleinbock's result.

\begin{theorem}[Zhang, 2009]\label{Z1}
If $\cL$ is an affine subspace of $\bbR^{n}$ and $\cM$ is a nondegenerate submanifold in $\cL$, then
\begin{equation}
\omega_{0}(\cM) = \omega_{0}(\cL) = \inf \{ \omega_{0}(\bx) \mid \bx \in \cL \} = \inf\{ \omega_{0}(\bx) \mid \bx \in \cM \} 
\end{equation}
Further if $\cL$ is a hyperplane defined by $\ba := (a_{1}, \ldots , a_{n-1}) \in \bbR^{n-1}$ as in \eqref{defL}, then
\begin{equation}
\omega_{0}(\cL)=\max \left\{ 1/n , \cfrac{\omega_{n-1}(\ba)}{n+(n-1)\omega_{n-1}(\ba)} \right\}.
\end{equation}
\end{theorem}

We should mention that some other cases of explicit computations of Diophantine exponents of subspaces have been calculated in \cite{Kleinbock-exponent} but that this problem is largely open and seems difficult. On the other hand, as far as we are aware, the corresponding inhomogeneous problem has not been studied so far and does not seem approachable directly using the techniques of Kleinbock and Zhang which are based on sharp nondivergence estimates for polynomial-like flows on the space of lattices developed by Kleinbock-Margulis and Kleinbock. 
In this paper, we extend Theorem \ref{BV} from extremal transfer to transfer for arbitrary exponents and as a consequence, obtain the first known bounds for the inhomogeneous Diophantine exponent of affine subspaces and their nondegenerate submanifolds.

\begin{theorem}\label{GM1u}
Let $\cM$ be a differentiable submanifold of $\bbR^{n}$. For every $\btheta\in\bbR^{n}$ we have that
\begin{equation}
\omega_{0}(\bx,\btheta) \leq \omega_{0}(\cM) \textrm{ for almost all } \bx \in \cM.
\end{equation}
 and
 \begin{equation}
\omega_{n-1}(\bx,\btheta) \leq \omega_{n-1}(\cM) \textrm{ for almost all } \bx \in \cM.
\end{equation}
\end{theorem}

In this case, a lower bound is still given by the transfer inequality of Bugeaud and Laurent \cite{BL}. It reads as follow.

\begin{theorem}\label{GM1l}
In the setting of Theorem \ref{GM1u} we also have

\begin{eqnarray}
\omega_{0}(\bx,\btheta) &\geq & \max( 0, 1-(n-1)\omega_{0}(\cM)) \textrm{ for all } \bx \in \cM.\\
\omega_{n-1}(\bx,\btheta) &\geq & \cfrac{\omega_{n-1}(\cM)}{\omega_{n-1}(\cM)-n+1} \textrm{ for all } \bx \in \cM.
\end{eqnarray}
\end{theorem}

We postpone the proof and discussion about these lower bounds to section \ref{Lowerbounds}. Here and later, we provide a lower bound for inhomogeneous exponents in term of their corresponding homogeneous exponent. This is not the case in the transfer results of Bugeaud and Laurent, so we combine their result with other transfer results due to German \cite{Ger1, Ger2}. We might lose optimality, but less information is required to apply our result. Remember that computing Diophantine exponents on an explicit example is a difficult problem.\\

Combining Theorems \ref{Z1}, \ref{KHypn-1}, \ref{GM1u} and \ref{GM1l}, we get the following corollary.

\begin{corollary}\label{application}
Let $\ba=(a_{1}, \ldots , a_{n-1})$ be a point in $\bbR^{n-1}$. Let $\cL$ be an hyperplane of $\bbR^{n}$ parametrized by $\ba$ as in (\eqref{defL}). Then, for any nondegenerate submanifold $\cM \subset \cL$, for every $\btheta\in\bbR^{n}$ and almost every $\bx \in \cM$ we have
\begin{eqnarray*}
  \min\left\{ 1/n, \cfrac{n}{n +(n-1)\omega_{n-1}(\ba)}\right\} &\leq& \omega_{0}(\bx,\btheta) \leq \max\left\{1/n , \cfrac{\omega_{n-1}(\ba)}{n +(n-1)\omega_{n-1}(\ba)}\right\}, \\
  \min\left\{ n , \cfrac{ns}{n \omega_{0}(\ba) -s(n-1)} \right\} &\leq& \omega_{n-1}(\bx,\btheta) \leq \max\left\{ n , \cfrac{\omega_{0}(\ba)}{\omega_{0}(\ba) -n+1}\right\}.
\end{eqnarray*}
\end{corollary}

In the next section, we present a more general, in particular, multiplicative setting, and provide in this context an extended version of our Theorem \ref{GM1u}.\\  

\section*{Acknowledgements}
We thank Y. Bugeaud, D. Kleinbock and S. Velani for helpful comments.

\section{A more general setting}

First, we need to define more general exponents of inhomogeneous Diophantine approximation.\\

Let $m,n \in \bbN$ and $\bbR^{m\times n}$ be the set of all $m \times n $ real matrices. Given $\bX \in \bbR^{m\times n}$ and $\btheta \in \bbR^{m}$, let $\omega(\bX,\btheta)$ be the supremum of $w\geq0$ such that for arbitrarily large $Q>1$ there exists a nonzero $\bq = (q_{1}, \ldots , q_{n})\in \bbZ^{n}$ satisfying
\begin{equation}\label{defexp}
\|\bX\bq + \btheta\| < Q^{-w} \textrm{ and } |\bq| \leq Q,
\end{equation}
where $|\bq| := \max\{ |q_{1}|, \ldots , |q_{n}| \}$ is the supremum norm and $\| \cdot \|$ is the distance to a nearest integer point. We denote by $\hat{\omega}(\bX,\btheta)$ the corresponding uniform exponent, that is the supremum of $w\geq0$ such that \eqref{defexp} has a solution for all $Q$ sufficiently large. Here and elsewhere, $\bq\in \bbZ^{n}$ and $\btheta\in \bbR^{m}$ are treated as columns. Note that we recover the exponents $\omega_{0}$ and $\omega_{n-1}$ when $m=1$ or $n=1$. Further, let us define the multiplicative exponents $\omega^{\times}(\bX, \btheta)$ (resp. $\hat{\omega}^{\times}(\bX, \btheta)$ ) to be the supremum of $w\geq0$ such that for arbitrarily large $Q>1$ (resp. every sufficiently large $Q$) there exists a nonzero $\bq := (q_{1}, \ldots , q_{n}) \in \bbZ^{n}$ satisfying
\begin{equation}\label{defexpmult}
\Pi \langle \bX \bq + \btheta \rangle < Q^{-mw} \textrm { and } \Pi_{+}(\bq) \leq Q^{n},
\end{equation}
where \[\Pi \by := \Pi (\by) = \prod_{j=1}^{m} |y_{j}| \textrm{ and } \Pi_{+}(\bq) := \prod_{i=1}^{n}\max\{ 1, |q_{i}|\} \]
for $\by =(y_{1}, \ldots , y_{m})$. Also, $\langle \by \rangle$ denotes the unique point in $[-1/2,1/2)^{m}$ congruent to $\by \in \bbR^{m}$ modulo $\bbZ^{m}$. Thus $\|\cdot\| = |\langle \cdot \rangle|$.\\

If $\btheta=0$ we are in the homogeneous setting. In this case, Dirichlet's pigeonhole principle provides that $\omega^{\times}(\bX) \geq \omega(\bX) \geq \tfrac{m}{n}$ for all $\bX\in\bbR^{m\times n}$. For Lebesgue almost all $\bX\in\bbR^{m\times n}$, the Borel-Cantelli lemma ensure that $\omega(\bX)=\tfrac{m}{n}$ and that $\omega^{\times}(\bX)=\tfrac{m}{n}$.\\

Note that Beresnevich and Velani use a different normalization, so that the 'extremal' value of each exponent is $1$. We chosed the normalization used in the transference principles from Bugeaud \& Laurent and German.\\

Subsequent to the work of Kleinbock and Margulis, a significant advance was made by Kleinbock, Lindenstrauss and Weiss \cite{KLW} where they defined the notion of ``friendly" measure and proved that almost every point in the support of such a measure is not very well multiplicatively approximable. The transference principles of Beresnevich and Velani are proved in the general context of (strongly) contracting measures, a category which includes friendly measures.\\ 

We follow the notation and terminology of Beresnevich and Velani \cite{BV1}. Let $\mu$ be a non-atomic, locally finite, Borel measure on $\bbR^{m+n}$. If $B$ is a ball in a metric space $\Omega$ then $cB$ denotes the ball with the same centre as $B$ and radius $c$ times the radius of $B$. A measure $\mu$ on $\Omega$ is non-atomic if the measure of any point in $\Omega$ is zero. The support of $\mu$ is the smallest closed set $S$ such $\mu(\Omega\backslash S) = 0$. Also, recall that $\mu$ is doubling if there is a constant $\lambda > 1$ such that for any ball $B$ with centre in $S$
$$\mu(2B) \leq\lambda \mu(B).$$
For $\ba \in \bbR^{n}$ with  $\|a\|_{2} =1$ and $\bb \in \bbR^{m}$ consider the plane
\begin{equation}\label{defplane}
\mathcal{L}_{\ba, \bb} := \{\bX \in \bbR^{m \times n} : \bX\ba + \bb = 0\}
\end{equation}

Given $\bepsilon = (\varepsilon_1,\dots, \varepsilon_m) \in (0, \infty)^m$, the $\epsilon$-neighborhood of the plane $\mathcal{L}_{\ba, \bb}$ is given by
\begin{equation}\label{defplanenbhd}
\mathcal{L}^{\bepsilon}_{\ba, \bb} := \{\bX \in \bbR^{m \times n} : |\bX_{j}\ba + \bb| < \varepsilon_j \text{ for all } 1 \leq j \leq m\},
\end{equation}
where $\bX_{j}$ is the $j$-th row of $\bX$. A non-atomic, finite, doubling Borel measure $\mu$ on $\bbR^{m\times n}$ is strongly contracting if there exist positive constants $C, \alpha$ and $r_0$ such that for any plane $\mathcal{L}_{\ba, \bb}$, any $\bepsilon = (\varepsilon_1,\dots, \varepsilon_m) \in (0, \infty)^m$ with $\min\{\varepsilon_j ~:~1 \leq j \leq m\} < r_0$ and any $\delta \in (0,1)$ the
following property is satisfied: for all $\bX \in \mathcal{L}^{\delta\bepsilon}_{\ba, \bb} \cap S$ there is an open ball $B$ centered at $\bX$ 
such that
\begin{equation}
B \cap S \subset \mathcal{L}^{\bepsilon}_{\ba, \bb}
\end{equation}
and
\begin{equation}
\mu(5B \cap \mathcal{L}^{\delta\bepsilon}_{\ba, \bb}) \leq C\delta^{\alpha}\mu(5B).
\end{equation}
The measure $\mu$ is said to be contracting if the property holds with $\varepsilon_1 = \dots = \varepsilon_m = \varepsilon$. We say that $\mu$ is (strongly) contracting almost everywhere if for $\mu$-almost every point $\bX_0 \in \bbR^{m\times n}$ there is a neighborhood $U$ of $\bX_0$ such that the restriction $\mu|_{U}$ of $\mu$ to $U$ is (strongly) contracting. The following Theorem is proved in \cite{BV1}. 

\begin{theorem}[Beresnevich, Velani, 2010]\label{BV1}
Let $\mu$ be a measure on $\bbR^{m \times n}$.
\begin{enumerate}
\item[(A)]If $\mu$ is contracting almost everywhere then
\begin{equation}\nonumber
\mu \text{ is extremal } \iff \mu \text{ is inhomogeneously extremal}.
\end{equation}
\item[(B)]If $\mu$ is strongly contracting almost everywhere then
\begin{equation}\nonumber
\mu \text{ is strongly extremal } \iff \mu \text{ is inhomogeneously strongly extremal}.
\end{equation}
\end{enumerate}
\end{theorem}

Similarly to the simpler case, the main content of \cite{BV1, BV2} is to provide an upper bound on the Diophantine exponent, and the corresponding lower bound is provided using a transference inequality of Bugeaud and Laurent \cite{BL}. We extend Theorem \ref{BV1} to the non-extremal case as follows.

\begin{theorem}\label{GM2u}
Let $\mu$ be a measure on $\bbR^{m \times n}$.
\begin{enumerate}
\item[(A)]If $\mu$ is contracting almost everywhere and if for $\mu$-almost every $\bX\in\bbR^{m \times n}$ we have $\omega(\bX)=v$ then for every $\btheta\in\bbR^{m}$ 
\begin{equation}\nonumber
\omega(\bX,\btheta) \leq v \textrm{ for $\mu$-almost every } \bX\in\bbR^{m \times n}
\end{equation}
\item[(B)]If $\mu$ is strongly contracting almost everywhere and if for $\mu$-almost every $\bX\in\bbR^{m \times n}$ we have $\omega^{\times}(\bX)=v^{\times}$ then for every $\btheta\in\bbR^{m}$ 
\begin{equation}\nonumber
\omega^{\times}(\bX,\btheta) \leq v^{\times} \textrm{ for $\mu$-almost every } \bX\in\bbR^{m \times n}.
\end{equation}
\end{enumerate}
\end{theorem}

In these settings, the lower bound reads as follows.

\begin{theorem}\label{GM2l}
Let $\mu$ be a measure on $\bbR^{m \times n}$.
\begin{enumerate}
\item[(A)]If $\mu$ is contracting almost everywhere and if for $\mu$-almost every $\bX\in\bbR^{m \times n}$ we have $\omega(\bX)=v$ then for every $\btheta\in\bbR^{m}$ and $\mu$-almost every $\bX\in\bbR^{m \times n}$
\begin{equation}\nonumber
\omega(\bX,\btheta) \geq 
\left\{\begin{array}{ll} \cfrac{v}{nv-m+1} & \textrm{ if } \hat{\omega}(^{t}\bX) \leq 1\\[4mm] m-(n-1)v& \textrm{ if } \hat{\omega}(^{t}\bX) \geq 1 \textrm{ e.g. } n\geq m\end{array}\right.
\end{equation}
\item[(B)]If $\mu$ is strongly contracting almost everywhere and if for $\mu$-almost every $\bX\in\bbR^{m \times n}$ we have $\omega^{\times}(\bX)=v^{\times}$ then for every $\btheta\in\bbR^{m}$ 
\begin{equation}\nonumber
\omega^{\times}(\bX,\btheta) \geq \max\left( 0, \cfrac{n-(m-1)v^{\times}}{mv^{\times} -n+1}\right)
\end{equation}
\end{enumerate}
\end{theorem}

These lower bounds are interesting whenever $v$ or $v^{\times}$ belong to the interval $\left[n/m, n/(m-1) \right]$.\\

Furthermore, Beresnevich and Velani show that any friendly measure on $\bbR^{n}$ is strongly contracting. Note that Riemannian measures supported on non-degenerate manifolds are known to be friendly \cite{KLW}.\\

Then, using a slicing argument, Beresnevich and Velani prove the following.

\begin{theorem}\label{BVslic}
Let $\cM$ be a differentiable submanifold of $\bbR^{n}$. Then 
\begin{enumerate}
\item[(A)]{Let $\cM$ be a differentiable submanifold of $\bbR^{n}$. Then 
\begin{equation}\nonumber
\cM \text{ is extremal } \iff \cM \text{ is simultaneously inhomogeneously extremal}.
\end{equation}}
\item[(B)]{Furthermore, suppose that at almost every point on $\cM$ the tangent plane is not orthogonal to any of the coordinate axes. Then
\begin{equation}\nonumber
\cM \text{ is strongly extremal } \iff \cM \text{ is simultaneously inhomogeneously strongly extremal}.
\end{equation}}
\end{enumerate}
\end{theorem}

Note that a measure supported on a differentiable manifold is not necessarily friendly. Also, (A) is in fact Theorem \ref{BV}, already extend to Theorems \ref{GM1u} and \ref{GM1l}. We extend the multiplicative result to the non-extremal case.

\begin{theorem}\label{GM3u}
Let $\cM$ be a differentiable submanifold of $\bbR^{n}$. Suppose that at almost every point on $\cM$ the tangent plane is not orthogonal to any of the coordinate axes. Then, for every $\btheta\in\bbR^{m}$ we have
\begin{eqnarray}\nonumber
\omega_{0}^{\times}(\bx,\btheta) \leq \omega_{0}^{\times}(\cM) \text{ for almost every } \bx \in \bbR^{n},\\
\omega_{n-1}^{\times}(\bx,\btheta) \leq \omega_{n-1}^{\times}(\cM) \text{ for almost every } \bx \in \bbR^{n}.
\end{eqnarray}
\end{theorem}

In this settings, the multiplicative lower bounds read as follow.

\begin{theorem}\label{GM3l}
With notation and conditions of Theorem \ref{GM3u}, we also have
\begin{eqnarray}\nonumber
\omega_{0}^{\times}(\bx,\btheta) \geq \cfrac{1-(n-1)\omega_{0}^{\times}(\cM)}{n\omega_{0}^{\times}(\cM)} \text{ for almost every } \bx \in \bbR^{n},\\
\omega_{n-1}^{\times}(\bx,\btheta) \geq \cfrac{n}{\omega_{n-1}^{\times}(\cM)-(n-1)} \text{ for almost every } \bx \in \bbR^{n}.
\end{eqnarray}
\end{theorem}

In the multiplicative setting, Zhang \cite{Zhang2} provides also an example of non-extremal manifolds.

\begin{theorem}[Zhang, 2010]\label{Z2}
If $\cL$ is a hyperplane of $\bbR^{n}$ and $\cM$ is a nondegenerate submanifold in $\cL$, then
\begin{equation}
\omega^{\times}_{n-1}(\cL) = \omega^{\times}_{n-1}(\cM) = \inf\left\{ \omega^{\times}_{}(\bx) \mid \bx\in \cL \right\} = \inf\left\{\omega^{\times}_{n-1}(\bx) \mid \bx \in \cM \right\}
\end{equation}
Furthermore, suppose that $\cL$ is defined by
\begin{equation}\label{defLmult}
(x_{1}, x_{2}, \ldots ,x_{n-1}) \mapsto (a_{1}x_{1} + a_{2}x_{2} + \cdots + a_{n-1}x_{n-1}+a_{n}, x_{1}, x_{2}, \ldots , x_{n-1})
\end{equation}
Denote $\ba:=(a_{1}, \ldots , a_{n})$ and suppose that $s-1$ is the number of nonzero elements in $\{a_{1}, \ldots , a_{n-1}\}$. Then we have
\begin{equation}
\omega^{\times}_{n-1}(\cL) = \max \left( n , \cfrac{n}{s}\omega_{0}(\ba) \right)
\end{equation}
\end{theorem}
Note that the two different definitions of hyperplane \eqref{defL} and \eqref{defLmult} are slightly different.\\

Combining Theorems \ref{GM2u}, \ref{GM2l} and \ref{Z2} we obtain the following multiplicative analogue of Corollary \ref{application}.

\begin{corollary}
If $\cL$ is a hyperplane of $\bbR^{n}$ defined by $\ba$ and \eqref{defLmult} and $\cM$ is a nondegenerate submanifold in $\cL$.  Suppose that at almost every point on $\cM$ the tangent plane is not orthogonal to any of the coordinate axes. Then for all $\btheta\in \bbR^{n}$ and almost every $\bx\in\cM$ we have
\begin{equation}
\min\left(n, \cfrac{ns}{n \omega_{0}(\ba) - (n-1)s}\right) \leq \omega_{n-1}^{\times}(\bx,\btheta) \leq \max\left(n,\cfrac{n}{s}\omega_{0}(\ba)\right).
\end{equation}
where $s-1$ is the number of nonzero numbers among the $n-1$ first coordinates of $\ba$.
\end{corollary}

Note that all the condition can be fulfilled only if $\dim( \cM) \leq s$.


\section{Lower bounds}\label{Lowerbounds}

In this section, we use different transference inequalities to provide the lower bounds of the Theorems \ref{GM1l}, \ref{GM2l} and \ref{GM3l}. These lower bounds essentially follow from a transference inequality of Bugeaud and Laurent \cite{BL}. We then use other transference inequality of German to express the lower bounds of the inhomogeneous exponents in term of their homogeneous analogues.\\ 

\begin{theorem}[Bugeaud, Laurent, 2005]\label{BL}
Let $\bx, \btheta \in \bbR^{n}$. Then 

\begin{equation}\label{BLeq}
\omega(\bX,\btheta) \geq \frac{1}{\hat{\omega}({}^{t}\bX)} \; \textrm{ and } \; \hat{\omega}(\bX,\btheta) \geq \frac{1}{\omega(\bX)}
\end{equation}
with equality in \ref{BLeq} for Lebesgue almost every $\btheta \in \bbR^{n}$.
\end{theorem}

For multiplicative exponents, we have the following consequence.

\begin{corollary}\label{BLmult}
Let $\bx, \btheta \in \bbR^{n}$. Then 
\begin{equation}\label{BLeq}
\omega^{\times}(\bX,\btheta) \geq \frac{1}{\hat{\omega}^{\times}({}^{t}\bX)} \; \textrm{ and } \; \hat{\omega}^{\times}(\bX,\btheta) \geq \frac{1}{\omega^{\times}({}^{t}\bX)}
\end{equation}
\end{corollary}

\noindent It comes from the fact that
\begin{equation}\label{}
\omega^{\times}(\bX,\btheta) \geq \omega(\bX,\btheta) \textrm{ for all } \bX \in \bbR^{m\times n} \textrm{ and all } \btheta \in \bbR^{m}.
\end{equation}

In the context of Theorem \ref{GM1l}, we use the following transference inequalities established by German \cite{Ger1}
\begin{theorem}\label{Ger1}
For every $\bx \in \mathbb{R}^{n}$, we have
\begin{equation}\label{GJ}
\cfrac{\hat{\omega}_{n-1}(\bx)-1}{(n-1)\hat{\omega}_{n-1}(\bx)} \leq \hat{\omega}_{0}(\bx) \leq \cfrac{\hat{\omega}_{n-1}(\bx) -(n-1)}{\hat{\omega}_{n-1}(\bx)}
\end{equation}

\end{theorem}

Combining it with Theorem \ref{BL}, we get that for every $\bx\in\bbR^{n}$ and every $\btheta\in\bbR^{n}$ we have 
\begin{eqnarray}
\omega_{n-1}(\bx,\btheta) &\geq& \cfrac{1}{\hat{\omega}_{0}(\bx)} \geq \cfrac{\hat{\omega}_{n-1}(\bx)}{\hat{\omega}_{n-1}(\bx) -n+1} \geq \cfrac{\omega_{n-1}(\bx)}{\omega_{n-1}(\bx) -n+1} \\
\omega_{0}(\bx,\btheta) &\geq& \cfrac{1}{\hat{\omega}_{n-1}(\bx)} \geq \max\left( 0, 1-(n-1)\omega_{0}(\bx)\right).
\end{eqnarray}
Note that the second inequality is non trivial if and only if $1/n \leq {\omega}_{0}(\bx) \leq 1/(n-1)$. This comes from the fact that Theorem \ref{Ger1} provides an upper constraint on $\hat{\omega}_{n-1}(\bx)$ in terms of $\hat{\omega}_{0}(\bx)$ if and only if $\hat{\omega}_{0}(\bx) \leq 1/(n-1)$. Fortunately, this fits well with our application to Theorem \ref{Z1}, because for any hyperplan $L$, the exponent $\omega(L)$ belongs to the range $[1/n,1/(n-1)]$.\\

Now consider the more general context of Theorem \ref{GM2l}, German's transference inequalities \cite{Ger1} read as follows.

\begin{theorem}[German, 2011]\label{Ger2}
For every $\bX\in\bbR^{m\times n}$, for every $\btheta\in\bbR^{m}$, we have
\begin{equation}
\hat{\omega}(^{t}\bX) \geq 
\left\{\begin{array}{ll} \cfrac{n-1}{m-\hat{\omega}(\bX)} & \textrm{ if } \hat{\omega}(\bX) \leq 1\\[4mm] \cfrac{n-(\hat{\omega}(\bX))^{-1}}{m-1}& \textrm{ if } \hat{\omega}(\bX) \geq 1\end{array}\right.
\end{equation}
\end{theorem}

Combining it with Theorem \ref{BL}, we get that for every $\bx\in\bbR^{m\times n}$ and every $\btheta\in\bbR^{m}$ we have

\begin{equation}
\omega(\bX,\btheta) \geq \cfrac{1}{\hat{\omega}(^{t}\bX)} \geq 
\left\{\begin{array}{ll} \cfrac{\omega(\bX)}{n\omega(\bX)-m+1} & \textrm{ if } \hat{\omega}(^{t}\bX) \leq 1\\[4mm] m-(n-1)\omega(\bX)& \textrm{ if } \hat{\omega}(^{t}\bX) \geq 1\end{array}\right.
\end{equation}

It is more interesting than Theorem \ref{BL} if we can get rid of the condition on $ \hat{\omega}(^{t}\bX)$. Namely, we have an interesting non trivial lower bound if $n\geq m$ and $m/n \leq \omega(\bX) \leq m/(n-1)$. Then,
\begin{equation}
\omega(\bX,\btheta) \geq m-(n-1)\omega(\bX) \geq 0
\end{equation}

In the multiplicative setting, we use an other set of transference inequalities stated by German \cite{Ger2}

\begin{theorem}[German, 2011]\label{Ger2}
Let $\bX \in \mathbb{R}^{m\times n}$, we have
\begin{equation}
\hat{\omega}^{\times}(\bX) \leq \cfrac{m\hat{\omega}^{\times}(^{t}\bX)-n+1}{n-(m-1)\hat{\omega}^{\times}(^{t}\bX)}.
\end{equation}
\end{theorem}

Combining it with Theorem \ref{BL}, we get for every $\bX\in\bbR^{m\times n}$ and every $\btheta\in\bbR^{m}$:
\begin{equation}
{\omega}^{\times}(\bX,\btheta)\geq \cfrac{1}{\hat{\omega}^{\times}(^{t}\bX)} \geq \max\left(0, \cfrac{m-(n-1){\omega}^{\times}(\bX)}{n{\omega}^{\times}(\bX)-m+1}\right)
\end{equation}

Again, this is non trivial if and only if $m/n \leq \hat{\omega}^{\times}(\bX) \leq m/(n-1)$. In particular, we have

\begin{eqnarray}
\omega_{n-1}^{\times}(\bx,\btheta) &\geq \cfrac{n}{\omega_{n-1}^{\times}(\bx) - (n-1) } ,\\
\omega_{0}^{\times}(\bx,\btheta) &\geq \cfrac{1-(n-1)\omega_{0}^{\times}(\bx)}{n\omega_{0}^{\times}(\bx)}.
\end{eqnarray}


\section{Proof of Theorems \ref{GM1u} and \ref{GM3u}}

We refer the reader to the proof of Theorem \ref{BVslic} in \cite[\S2.3]{BV1}. Here, we just give a sketch and explain how to adapt it to the non-extremal case.\\

The idea is to apply Theorem \ref{GM2u}. Given a differential submanifold $\cM$ of $\bbR^{n}$, if we denote by $m$ the Riemannian measure on $\cM$, we only need to prove that $m$ is strongly contracting almost everywhere. We reduce the problem to the case of curves with a slicing argument. Once the result proved for the curves, we use Fubini's theorem to recover it for the whole manifold $\cM$.\\ 

Every step of the proof are the same as in \cite[\S2.3]{BV1}, we refer the reader to it. To adapt it to the non-extremal case, we just need to replace the set of full measure $\cE$ by either 
\[ \cE := \left\{ \bx \in B_{0} : \omega_{0}^{\times}(\bx)=\omega_{0}^{\times}(\cM)\right\} \textrm{ or } \cE := \left\{ \bx \in B_{0} : \omega_{n-1}^{\times}(\bx)=\omega_{n-1}^{\times}(\cM)\right\}.\]
and at the end with Fubini's theorem we prove that either 
\[ \cE^{\btheta} := \left\{ \bx \in B_{0} : \omega_{0}^{\times}(\f(\bx),\btheta)=\omega_{0}^{\times}(\cM)\right\} \textrm{ or } \cE^{\btheta} := \left\{ \bx \in B_{0} : \omega_{n-1}^{\times}(\f(\bx),\btheta)=\omega_{n-1}^{\times}(\cM)\right\}\]
has full dimension.\\

It is also possible to get a self-contained proof of Theorem \ref{GM1u} by adapting the proof from \cite{BV2} in a similar way.\\


\section{Proof of Theorem \ref{GM2u}}

\subsection{A reformulation of Theorem \ref{GM2u}}
Following the steps of \cite{BV1}, we introduce some notations adapted to the non-extremal setting and reformulate Theorem \ref{GM2u}. Then, we state the transference theorem of Beresnevich and Velani in its full bright and use it for our proof.\\

Let $\mu$ be a strongly extremal measure on $\bbR^{m\times n}$ and define the set

\[ \cA_{m,n}^{\btheta}(v^{\times}) := \left\{ \bX \in \bbR^{m\times n} : \omega^{\times}(\bX,\btheta) > v^{\times} \right\}. \]

We prove Theorem \ref{GM2u} if we show that
\begin{equation}
\mu(\cA_{m,n}^{\btheta}(v^{\times}))=0 \quad \textrm{ for all } \btheta\in \bbR^{m}
\end{equation}

Let $\bT$ denote a countable subset of
$\bbR^{m+n}$ such that for every $\ttt=(t_1,\dots,t_{m+n})\in\bT$

\begin{equation}\label{condt1}
    \sum_{j=1}^m t_j=\lambda\sum_{i=1}^{n} t_{m+i}\,.
\end{equation}
where $\lambda := \tfrac{m}{n}v^{\times}$.\\

\noindent For $\ttt\in\bT$, consider the diagonal transformation $g_{\ttt} $ of $\bbR^{m+n}$ given by
\begin{equation}\label{e:040}
  g_{\ttt} \ := \ \operatorname{diag}\{2^{t_1},\dots,2^{t_m},2^{-t_{m+1}},\dots,2^{-t_{m+n}}\}\,.
\end{equation}

\noindent For $\bX \in \bbR^{m \times n}$, define the matrix
\[ M_{\bX}:=\left(\begin{array}{cc}
  I_m&\bX\\[2ex]
  0&I_n
\end{array}
 \right)\,,\]
where $I_n$ and $I_m$ are respectively the  $n\times n$ and $m\times m$ identity matrices. The matrix $M_{\bX}$ is a linear transformation of $\bbR^{m+n}$. Given $\btheta \in \bbR^m$, let
\[ M_{\bX}^{\btheta}\ :\ \ba\mapsto M_{\bX}^{\btheta}\ba:=M_{\bX}\ba+\bm\Theta\,,\]

\noindent where $\bm\Theta \, := \,
{}^t(\theta_1,\dots,\theta_m,0,\dots,0)\in\bbR^{m+n}$. Thus, $M_{\bX}^{\btheta}$ is  an affine transformation of $\bbR^{m+n}$.

\noindent Let
\begin{equation}\label{DefcA}
\cA=\bbZ^m\times(\bbZ^n\setminus \{ \bf{0} \})\,.
\end{equation}
Then, for  $\varepsilon>0$, $\ttt\in\bT$ and $\alpha\in\cA$ define the sets
\begin{equation}\label{DefDelta}
\Delta^{\btheta}_\ttt(\alpha,\varepsilon):=\{\bX\in\bbR^{m \times n}:| g_{\bt}M_{\bX}^{\btheta}\alpha|\,< \varepsilon\}
\end{equation}
and
$$
\Delta^{\btheta}_\ttt(\varepsilon):=\bigcup_{\alpha\in
\cA}\Delta^{\btheta}_\ttt(\alpha,\varepsilon)= \{\bX\in\bbR^{m \times n}:\inf_{\alpha\in\cA}|
g_{\ttt}M_{\bX}^{\btheta}\alpha|\,< \varepsilon\}\,.
$$

\noindent For $\eta>0$,  define the function
\begin{equation}\label{DefPsi}
\psi^\eta \, : \, \bT \mapsto \bbR_+   \  :  \
\ttt\mapsto\psi^\eta_\ttt :=2^{-\eta\sigma(\ttt)} \
\end{equation}
where $\sigma(\ttt):=t_1+\dots+t_{m+n}$, and consider the $\limsup$ set given by
\begin{equation}\label{DefLa}
 \La^{\btheta}_{\bT}(\psi^\eta\,) \, := \, \limsup_{\ttt \in \bT }\Delta^{\btheta}_\ttt(\psi^\eta_\ttt)\, .
\end{equation}
In the case $\btheta=\bf{0}$, we write $\La_{\bT}(\psi^\eta)$ for $\La^{\btheta}_{\bT}(\psi^\eta)$. The following result provides a reformulation of the set $ \cA_{m,n}^{\btheta} $ in terms of the $\limsup$ sets given by  \eqref{DefLa}.

\begin{proposition}\label{Key}
There exists a countable subset $\bT$ of $\bbR^{m\times n}$ satisfying \eqref{condt1} such that
\begin{equation}\label{condt2}
\sum_{\ttt\in\bT} 2^{-\eta\sigma(\ttt)} < \infty \qquad \forall \eta >0
\end{equation}
and
\begin{equation}\label{condt3}
\cA_{m,n}^{\btheta} (v^{\times})= \bigcup_{\eta>0} \La_{\bT}^{\btheta}( \psi^{\eta}) \qquad \forall \btheta \in \bbR^{m}
\end{equation}
\end{proposition}

In fact, in the proof of Proposition \ref{Key} we show that we can construct a set $\bT$ that fits the non-extremal setting \eqref{condt1} but still has the properties \eqref{condt2} and \eqref{condt3}. This is the key point in the extension of Theorem \ref{BVslic} to the non-extremal case. Thereafter, our limsup sets have the necessary properties to apply the Inhomogeneous Transference Principle. Namely, we are reduced to show that for a set $\bT$ given by Proposition \ref{Key},

\begin{equation}\label{equivthm}
\mu( \La_{\bT}(\psi^{\eta}) ) = 0 \quad \forall\eta>0 \quad \Longrightarrow \quad  \mu( \La_{\bT}^{\btheta}(\psi^{\eta}) ) = 0 \quad \forall\eta>0
\end{equation} 

\subsection{Proof of Proposition \ref{Key}}

Given $\bs=(s_{1}, \ldots , s_{m}) \in \bbZ^{m}_{+}$ and $\bll=(l_{1}, \ldots , l_{n})\in \bbZ^{n}_{+}$, let 
\begin{equation*}
\sigma(\bs) := \sum_{j=1}^{m}s_{j} , \quad \sigma(\bll) := \sum_{i=1}^{n}l_{i} \quad \textrm{and} \quad \zeta:=\zeta(\bs,\bll) = \cfrac{\sigma(\bs)-\lambda \sigma(\bll)}{m+\lambda n},
\end{equation*}
where $\bbZ_{+}$ is the set of non-negative integers. Furthermore,
define the $(m+n)$-tuple $\ttt=(t_1,\dots,t_{m+n})$ by setting
\begin{equation}\label{deft}
\ttt := \left(s_{1}- \zeta , \ldots , s_{m} - \zeta, l_{1}+\zeta , \ldots , l_{n}+\zeta \right)
\end{equation}
and let
\begin{equation}\label{defT}
\bT := \left\{\ttt \in \bbR^{m\times n} \textrm{ defined by \eqref{deft} } : \bs \in \bbZ^{m}_{+}, \bll \in \bbZ^{n}_{+} \textrm{ with } \sigma(\bs) \geq \lambda\sigma(\bl)\right\}.
\end{equation}
We aim at showing that this choice of $\bT$ is suitable within the context of Proposition \ref{Key}. The choice of $\zeta$ ensure that definition \eqref{deft} satisfies condition \eqref{condt1}. First, we check that this $\bT$ satisfies \eqref{condt2}. For any $\ttt \in \bT$,
\begin{equation}\label{calcsigma12}
\tfrac{\lambda}{1+\lambda} \sigma(\ttt) = \sigma(\bs) - m \zeta \quad \textrm{ and } \quad \tfrac{1}{\lambda +1} \sigma(\ttt) = \sigma(\bll)+n\zeta
\end{equation}
where $\sigma(\ttt) := \sum_{k=1}^{m+n}t_{k}$. Since $\zeta$ is non-negative, we deduce that
\begin{equation}\label{relsigma}
(\lambda+1) \sigma(\bll) \leq \sigma(\ttt) \leq  \tfrac{\lambda+1}{\lambda}\sigma(\bs).
\end{equation}
Furthermore, on summing the two expressions arising in \eqref{calcsigma12} and using the fact that $\sigma(\bll)\geq0$, we obtain that
\begin{equation}\label{minsigmat}
\begin{array}[b]{rcl}
\sigma(\ttt)  = \sigma(\bs) +\sigma(\bl)+(m-n)\zeta   & \geq & \cfrac{\lambda+1}{m+n\lambda}\left(\sigma(\bs)+ \sigma(\bll) \right) .
\end{array}
\end{equation}
This ensures that $\bT$ satisfy condition \eqref{condt2}. In turn, it follows that for any $v\in \bbR_{+}$
\begin{equation}
\# \{ \ttt \in \bT : \sigma(\ttt)< v \} < \infty
\end{equation}

Now we check condition \eqref{condt3}. Fix $\btheta \in \bbR^{m}$. Note that $\bX \in \cA_{m,n}^{\btheta}(v^{\times})$ if and only if there exists $\varepsilon>0$, such that for arbitrarily large $Q>1$ there is an $\alpha=(\bp,\bq)\in \cA := \bbZ^{m} \times( \bbZ^{n}\setminus \{\bf{0}\})$ satisfying $|\bX\bq + \bp +\btheta| \leq 1/2$ such that
\begin{equation}\label{condens}
\Pi(\bX\bq+\bp+\btheta) < Q^{(1+\varepsilon)mv^{\times}} \quad \textrm{ and } \qquad \Pi_{+}(\bq) \leq Q^{n}.
\end{equation}

\paragraph{\textbf{Step 1.}} We show the inclusion
\begin{equation}\label{incl1}
\cA_{m,n}^{\btheta}(v^{\times}) \subseteq \bigcup_{\eta>0}\La_{\bT}^{\btheta}(\psi^{\eta}). 
\end{equation}
Suppose $\bX \in \cA_{m,n}^{\btheta}(v^{\times})$. It follows that \eqref{condens} is satisfied fo infinitely many $Q \in \bbZ_{+}$. For each such $Q$, we consider the unique $\bs\in\bbZ^{m}_{+}$ and $\bll\in\bbZ_{+}^{n}$ such that
\begin{equation}\label{defsl1}
    2^{-s_j}\leq \max\Big\{|\bX_{j} \, \bq+p_j+\theta_j|\ ,\
    Q^{-(1+\varepsilon)}\Big\}<2^{-s_j+1}\qquad\text{for } \ \ 1\leq j\leq m
\end{equation}
and
\begin{equation}\label{defsl2}
    2^{\, l_i}\leq \max\{1,|q_i|\}<2^{\, l_i+1}\qquad\text{for } \ \ 1\leq i\leq
    n\,.
\end{equation}
Here and after, $\bX_{j}:=(x_{j,1},\dots,x_{j,n})$ denotes the $j$-th row of $\bX\in\bbR^{m\times n}$. If we multiply over the indexes we get 
\begin{eqnarray}
2^{\sigma(\bll)} &\leq& \Pi_{+}(\bq) \leq Q^{n}, \\
2^{-\sigma(\bs)} &<& \max\left\{\Pi(\bX\bq+\bp+\btheta)   ,Q^{-(1+\varepsilon)mv^{\times}}
\right\} = Q^{-(1+\varepsilon)mv^{\times}}
\end{eqnarray}
\noindent Combining both inequalities, we get $2^{-\sigma(\bs)} < 2^{\sigma(\bll)(1+\varepsilon)mv^{\times}}$. Hence,
\begin{equation}
\sigma(\bs) - \lambda \sigma(\bll) > \varepsilon \lambda \sigma(\bll) \geq 0.
\end{equation}
Thus, $\ttt$ given by \eqref{deft} with $\bs$ and $\bll$ as defined above in \eqref{defsl1} and \eqref{defsl2} belongs to $\bT$.\\

\noindent If $\sigma(\bs) > 2 \lambda \sigma(\bll)$, then
\[ \zeta = \cfrac{\sigma(\bs)-\lambda\sigma(\bll)}{m+n\lambda} \geq \cfrac{\sigma(\bs)}{2(m+n \lambda)}\stackrel{\eqref{relsigma}}{\geq} \cfrac{\lambda\sigma(\ttt)}{2(\lambda+1)(m+n\lambda)}.\]
If $\sigma(\bs) \leq 2 \lambda\sigma(\bll)$, then
\[ \zeta = \cfrac{\sigma(\bs)-\lambda\sigma(\bll)}{m+n\lambda} \geq  \cfrac{\varepsilon\lambda\sigma(\bll)}{m + n \lambda} \geq \cfrac{\varepsilon \sigma(\bs)}{2(m+n\lambda)}\stackrel{\eqref{relsigma}}{\geq} \cfrac{\varepsilon\lambda\sigma(\ttt)}{2(\lambda+1)(m+n\lambda)} .\]
On combining the two cases, we deduce that
\begin{equation}\label{lowerzeta}
\zeta > \eta_{0} \sigma(\ttt) \qquad \textrm{ with } \quad \eta_{0}:= \cfrac{ \lambda}{2(\lambda+1)(m+n\lambda)}\min\left( 1,\varepsilon \right)
\end{equation}
The diagonal transformation $g_{\ttt}$ satisfies
\[ g_{\ttt} = 2^{-\zeta} \operatorname{diag}\{2^{s_{1}}, \ldots , 2^{s_{m}}, 2^{-l_{1}}, \ldots , 2^{-l_{n}} \}\]
It follow from definition \eqref{defsl1} and \eqref{defsl2} that
\begin{equation}
\inf_{\alpha\in\cA} | g_{\ttt} M_{\bX}^{\btheta}\alpha | < 2\cdot 2^{-\zeta}.
\end{equation}
For $0 < \eta \leq \eta_{0}$, the lower bound \eqref{lowerzeta} for $\zeta$ implies that
\begin{equation}\label{majgMa}
\inf_{\alpha\in\cA} | g_{\ttt} M_{\bX}^{\btheta}\alpha | < 2^{-\eta\sigma(\ttt)}
\end{equation}
for all sufficiently large $\sigma(\bt)$. Note that \eqref{condens} and \eqref{defsl1} ensure that $\sigma(\bs) \to \infty$ as $Q \to \infty$. 
Since \eqref{condens} is satisfied for arbitrarily large $Q\in\bbZ_{+}$ and \eqref{minsigmat} ensures that $\sigma(\ttt)$ also goes to infinity with $Q$, we have that \eqref{majgMa} is satisfied for infinitely many $\ttt \in \bT$. This proves that $\bX \in \La_{\bT}^{\btheta}(\psi^{\eta})$ for any $\eta \in (0,\eta_{0})$. This establishes the inclusion \eqref{incl1}.\\

\paragraph{Step 2} We show the inclusion
\begin{equation}\label{incl2}
\cA_{m,n}^{\btheta}(v^{\times})\supseteq \bigcup_{\eta>0}\La_{\bT}^{\btheta}(\psi^{\eta}). 
\end{equation}

Suppose that $\bX\in \La_{\bT}^{\btheta}(\psi^{\eta})$ for some $\eta>0$. By definition, \eqref{majgMa} is satisfied for infinitely many $\ttt \in \bT$. For each such $\ttt$, there exists $\alpha=(\bp,\bq)\in \cA$ such that 
\begin{equation}
| g_{\ttt} M_{\bX}^{\btheta}\alpha | < 2^{-\eta\sigma(\ttt)}
\end{equation}
If we take the product over the first $m$ coordinates of $ g_{\ttt} M_{\bX}^{\btheta}\alpha $, we obtain that
\[ \prod_{j=1}^{m}2^{t_{j}}|\bX_{j}\bq + p_{j}+ \theta_{j}| < 2^{-m\eta\sigma(\ttt)}.\]
Similarily, the product of the last $n$ non-zero coordinates of $ g_{\ttt} M_{\bX}^{\btheta}\alpha $ gives that
\[\prod_{ \underset{q_{i}\neq0}{1\leq i \leq n}} 2^{-t_{m+i}}|q_{i}| < 2^{-n \eta \sigma(\ttt)}.\]
By definition, for every $\ttt\in\bT$, we have $t_{m+i} \geq 0$ ($1 \leq i \leq n$). Also, the minoration \eqref{minsigmat} ensure that $\sigma(\ttt) \geq 0$. We obtain
\begin{equation}
\Pi(\bX\bq +\bp + \btheta) < 2^{-m\eta\sigma(\ttt) - \frac{\lambda\sigma(\ttt)}{1+ \lambda}} \quad \textrm{ and } \quad \Pi_{+}(\bq) < 2^{-n\eta\sigma(\ttt) + \frac{\sigma(\ttt)}{1+\lambda}}
\end{equation}
If we set 
\[ Q := 2^{\frac{\sigma(\ttt)}{n(1+\lambda)}} \quad \textrm{ and } \quad \varepsilon := \cfrac{\lambda+1}{\lambda}m\eta, \]
it follows that \eqref{condens} is satisfied for arbitrarily large $Q$ and arbitrarily small $\varepsilon$. Hence, $\bX \in \cA_{m,n}^{\btheta}(v^{\times})$. This establishes \eqref{incl2}.\\

Steps $1$ and $2$ establish \eqref{condt3} and complete the proof of Proposition \ref{Key}.\qed\\

\emph{Remark:} With $\cA_{m,n}^{\btheta}(v) := \left\{ \bX \in \bbR^{m \times n}: \omega(\bX,\btheta) > v\right\}$, which is the setting for Theorem \ref{GM1u}, the proof is essentially the same. We just add the extra condition that $s_{1}= \cdots = s_{m}$ and $l_{1}= \cdots = l_{n}$ in the definition of $\bT$. As a subset of the previous set, it satisfies conditions \eqref{condt1} and \eqref{condt2}. Replacing \eqref{condens} by 
\[ \| \bX \bq + \bp + \btheta\| < Q^{-(1+\varepsilon) v} \quad \textrm{ and }\quad |\bq|<Q,\]
the arguments of Steps $1$ and $2$ can naturally be modified to obtain \eqref{condt3}.\\


\subsection{An Inhomogeneous Transference Principle}

We recall here the general framework of the transference theorem of Beresnevich and Velani as it appears in \cite{BV1}. It allows to transfert zero mesure statement for homogeneous limsup sets to inhomogeneous limsup sets.\\

Let $(\Omega,d)$ be a locally compact metric space. Given two countable `indexing' sets $\cA$ and $\bT$, let $H$ and $I$ be two maps from $\bT\times \cA\times\bbR^+ $ into the set of open subsets of $\Omega$ such that
\[ H\,:\,(\ttt,\alpha,\varepsilon)\in \bT\times \cA\times\bbR^+ \,\mapsto
\,H_\ttt(\alpha,\varepsilon)\]
and
\[ I\,:\,(\ttt,\alpha,\varepsilon)\in \bT\times \cA\times\bbR^+ \,\mapsto
\,I_\ttt(\alpha,\varepsilon)\,.\]

\noindent Furthermore, let
\begin{equation}
H_\ttt(\varepsilon):=\bigcup_{\alpha\in \cA}H_\ttt(\alpha,\varepsilon)\qquad \textrm{ and } \qquad
I_\ttt(\varepsilon):=\bigcup_{\alpha\in \cA}I_\ttt(\alpha,\varepsilon)\,.
\end{equation}

\noindent Next, let $\bm\Psi$ denote a set of functions $
\psi:\bT\to\bbR^+\,:\,\ttt\mapsto \psi_\ttt\,. $ For $\psi\in\bm\Psi$,
consider the $\limsup$ sets

\begin{equation}
 \La_H(\psi\,)=\limsup_{\ttt \in \bT}H_\ttt(\psi_\ttt)
 \qquad \textrm{ and } \qquad
 \La_I(\psi\,)=\limsup_{\ttt \in \bT}I_\ttt(\psi_\ttt)\,.
\end{equation}

\bigskip

\noindent For reasons that will soon become apparent, we refer to sets associated with the map $H$  as homogeneous sets and those
associated with the map $I$ as inhomogeneous sets. The following `intersection' property states that the intersection of two distinct inhomogeneous sets is contained in a homogeneous set.

\bigskip

\noindent\textbf{The intersection property. } The triple $(H,I,\bm\Psi)$ is said to satisfy \emph{the intersection property} if for any $\psi\in\bm\Psi$, there exists $\psi^*\in\bm\Psi$ such that for all but finitely many $\ttt\in\bT$ and all  distinct $\alpha$ and $\alpha'$ in $\cA$ we have that 
\begin{equation}
    I_\ttt(\alpha,\psi_\ttt)\cap I_\ttt(\alpha',\psi_\ttt)\subset
H_{\ttt}(\psi^*_{\ttt})   \ .
\end{equation}

\bigskip

\noindent\textbf{The contracting property. } Let $\mu $ be a non-atomic, finite, doubling  measure supported on a bounded subset
$\bS$ of $\Omega$.  We say that  $\mu$ is \emph{contracting with respect to\/ $(\,I,\bm\Psi)$}\/ if for any $\psi\in\bm\Psi$ there
exists $\psi^+\in\bm\Psi$ and a sequence of positive numbers $\{k_\ttt\}_{\ttt\in \bT}$ satisfying
\begin{equation}
    \sum_{\ttt\in \bT}k_\ttt<\infty  \ ,
\end{equation}
 such that for all but finitely  $\ttt\in \bT$ and all $\alpha\in \cA$ there exists a collection $\mathcal{C}_{\ttt,\alpha}$ of balls $B$
centred at $\bS$ satisfying the following conditions\,{\rm:}
  \begin{equation}
    \bS\cap I_\ttt(\alpha,\psi_\ttt) \ \subset \
    \bigcup_{B\in\mathcal{C}_{\ttt,\alpha}}B\
  \end{equation}
    \begin{equation}
        \bS\cap\bigcup_{B\in\mathcal{C}_{\ttt,\alpha}}B \ \subset \ I_\ttt(\alpha,\psi^+_{\ttt})
    \end{equation}
    and
    \begin{equation}
        \mu\Big(5B\cap I_\ttt(\alpha,\psi_\ttt)\Big)\ \le  \ k_\ttt\ \,  \mu(5B) \ .
    \end{equation}

\vspace{4ex}

The intersection and contracting properties enable us to transfer zero $\mu$-measure statements for the homogeneous $\limsup$ sets
$\La_H(\psi\,)$ to the inhomogeneous $\limsup$ sets $\La_I(\psi\,) $.

\begin{theorem}[Inhomogeneous Transference Principle]\label{BVtransfert} 
Suppose that $(H,I,\boldsymbol{\Psi})$ satisfies the intersection property and that $\mu$ is contracting with respect to $(I,\boldsymbol{\Psi})$. Then
\begin{equation}
\mu(\La_{H}(\psi)) = 0 \quad \forall \psi \in \boldsymbol{\Psi} \Rightarrow \mu(\La_{I}(\psi)) = 0 \quad \forall \psi \in \boldsymbol{\Psi}
\end{equation}
\end{theorem}

\bigskip 

\subsection{Conclusion of the proof}

Throughout $\btheta\in\bbR^m$ is fixed.  Let $\mu$ be a measure on $\bbR^{m \times n}$ that is strongly contracting almost everywhere and fix a set $\bT$ arising from Proposition \ref{Key}. In terms of establishing \eqref{equivthm}, sets of $\mu$-measure zero are irrelevant. Therefore we can simply assume that $\mu$ is strongly contracting. We  show that \eqref{equivthm} falls within the scope of the above general framework. Let $\Omega := \bbR^{m \times n}$ and let $\cA$ be given by \eqref{DefcA}. Given $\varepsilon \in \bbR^{+}$,  $\ttt \in \bT$ and $ \alpha \in \cA$  let
\[H_\ttt(\alpha,\varepsilon):=\Delta_\ttt(\alpha,\varepsilon)=\Delta^{0}_\ttt(\alpha,\varepsilon)\qquad\text{and}\qquad
 I_\ttt(\alpha,\varepsilon):=\Delta^{\btheta}_\ttt(\alpha,\varepsilon),\]
where $\Delta^{\btheta}_\ttt(\alpha,\varepsilon)$ is defined by \eqref{DefDelta}.
This defines the  maps $H$ and $I$ associated with the general framework.  It is readily seen that
$H_\ttt(\varepsilon)=\Delta^{0}_\ttt(\varepsilon)$ and
$I_\ttt(\varepsilon)=\Delta^{\btheta}_\ttt(\varepsilon)$. Next, let $\bm\Psi$ be the class of  functions given by \eqref{DefPsi}. Then,  it immediately follows that
\[  \La_H(\psi) \,  = \, \La_{\bT}(\psi) := \La^{0}_{\bT}(\psi)
  \qquad\text{and}\qquad \La_I(\psi) \, = \,   \La^{\btheta}_{\bT}(\psi) \ ,
\]
where the set $\La^{\btheta}_{\bT}(\psi)$ is  defined by \eqref{DefLa}. In \cite{BV1}, it is shown that these sets satisfy the \emph{intersection property} and the \emph{contracting property}. Namely, only \eqref{condt1} changes in the non-extremal setting, and it is used only used once to show that $\sigma(\ttt)\geq 0$ implies $\sum_{j=1}^{m}t_{j} = \tfrac{\lambda}{1+ \lambda}\sigma(\ttt) \geq 0$. This remains true if $\lambda >1$. Thus we can apply Theorem \ref{BVtransfert}, which proves Theorem \ref{GM2u}. \qed\\

\section{Open problems}

The authors would like to point out that Beresnevich and Velani finish their paper \cite[\S8]{BV1} with a long and interesting presentation of open questions related to their Inhomogeneous Transference Principle, and strongly encourage the reader to look at it.\\

Concerning the non-extremal case, it would be interesting to provide explicite examples of manifolds where homogeneous and inhomogeneous exponents can be computed in order to check wether the inequality in Theorems \ref{GM1u} and \ref{GM1l} and Theorems \ref{GM3u} and \ref{GM3l} are best possible, and if the whole intervals is reachable.\\

In another direction, the theory of Diophantine approximation on manifolds discussed so far can be generalized to the context of smooth submanifolds of matrices, i.e. one considers submanifolds of systems of linear forms. The present theory corresponds to the special case of $n \times 1$ matrices. We refer the reader to \cite{KMW, BKM, ABRS, DFSU} for recent developments on this theme. One of the main difficulties in studying Diophantine approximation on submanifolds of matrices is that it doesn't seem straightforward to define the correct notion of nondegeneracy for submanifolds or indeed the right generalization of friendly measures. Accordingly, in the papers mentioned above, several notions have been developed to address this issue - for instance in \cite{BKM}, Beresnevich, Kleinbock and Margulis develop a notion of ``weakly non-planar" measures and in addition to proving the analogue of the Baker-Sprind\v{z}uk conjectures for such measures, an inhomogeneous transference principle is also proved, for the critical exponent, thereby generalising the work of Beresnevich and Velani. The results in the present paper also extend to this setting, however, we have chosen to restrict ourselves to the setting of measures and submanifolds of $\bbR^n$ because our results are especially significant for affine subspaces and the corresponding theory in the matrix setting is not yet sufficiently well developed even in the homogeneous approximation case.

Finally, to properly use Theorems \ref{GM2u} and \ref{GM2l} , it would be good to look for measures that are (strongly) contracting but not friendly.\\


\bibliographystyle{amsplain}

\end{document}